\documentclass[thmsa,a4paper,11pt]{article}%
\usepackage{amssymb}
\usepackage{amsmath}%
\usepackage{amsfonts}%
\usepackage{graphicx}
\providecommand{\U}[1]{\protect\rule{.1in}{.1in}}
\begin{document}

\author{Dirk Veestraeten\thanks{Address correspondence to Dirk Veestraeten, Amsterdam
School of Economics, Roetersstraat 11, 1018 WB Amsterdam, The Netherlands,
e-mail: d.j.m.veestraeten@uva.nl.}\\University of Amsterdam\\Amsterdam School of Economics\\Roetersstraat 11\\1018 WB Amsterdam\\The Netherlands\\\bigskip\bigskip}
\title{Integral representations for products of two parabolic cylinder functions with
different arguments and orders\\\bigskip\bigskip}
\maketitle

\begin{abstract}
\noindent This paper derives new integral representations for products of two
parabolic cylinder functions. In particular, expressions are obtained for
$D_{\nu}\left(  x\right)  D_{\mu}\left(  y\right)  $, with $x\geqslant0$ and
$y\geqslant0$, that allow for different orders and arguments in the two
parabolic cylinder functions. Also, two integral representations are obtained
for $D_{\nu}\left(  -x\right)  D_{\mu}\left(  y\right)  $ by employing the
connection between the parabolic cylinder function and the Kummer confluent
hypergeometric function. The integral representations are specialized for
products of two complementary error functions and of two modified Bessel
functions of the second kind of order $1/4$, as well as for the product of a
parabolic cylinder function and a modified Bessel function of the first kind
of order $1/4$.

\bigskip\bigskip

\bigskip


\end{abstract}

\newpage\setlength{\baselineskip}{1.4\baselineskip}

\section{Introduction}

The study of Nicholson--type integrals for the product of two parabolic
cylinder functions recently gained renewed interest. An integral
representation for the product of two parabolic cylinder functions with the
same order and identical or opposite arguments, $D_{\nu}(x)D_{\nu}(\pm x)$,
was presented in \cite{m03}. Recently, \cite{g15} obtained an integral
representation for the case of identical orders but unrelated arguments,
$D_{\nu}(x)D_{\nu}(y)$. Integral representations, on the other hand, for
unrelated orders but identical or opposite arguments, $D_{\nu}(\pm
x)D_{\nu+\mu-1}(x)$, were documented in \cite{n15}. The approach in
\cite{v15a} offered a first, yet incomplete attempt at allowing both for
differing arguments as well as orders given that the orders in the various
expressions still had to be linearly related. In particular, \cite{v15a} used
the results of \cite{v15b} to first obtain integral representations for
$D_{\nu}(x)D_{\nu}(y)$ and $D_{\nu}(x)D_{\nu-1}(y)$ that subsequently can be
extended via the recurrence relation for the parabolic cylinder function into
separate integral representations for other, linearly related orders such as
$D_{\nu}(x)D_{\nu+1}(y)$ and $D_{\nu}(x)D_{\nu-2}(y)$. This paper generalizes
existing results by deriving an integral representation for $D_{\nu}(x)D_{\mu
}(y)$ in which both the arguments as well as the orders are (linearly) unrelated.

The integral representation for $D_{\nu}\left(  x\right)  D_{\mu}\left(
y\right)  $, with $x\geqslant0$ and $y\geqslant0$, is derived by applying the
convolution theorem to two inverse Laplace transforms for single parabolic
cylinder functions that are documented in \cite{emoti154}. The latter
expressions are selected upon their ability to ultimately yield a Gaussian
hypergeometric function within the integrand of the integral representation.
The properties of the Gaussian hypergeometric function alternatively also
allow expressing the integrand of $D_{\nu}\left(  x\right)  D_{\mu}\left(
y\right)  $ in terms of the associated Legendre function of the first kind.
Subsequently, two expressions are derived for $D_{\nu}\left(  -x\right)
D_{\mu}\left(  y\right)  $ by employing the definition of the parabolic
cylinder function as a sum of two Kummer confluent hypergeometric functions.

The results of this paper can be specialized in various directions. For
instance, integral representations for products of two complementary error
functions and of two modified Bessel functions of the second kind of order
$1/4$ can immediately be obtained from the properties of the parabolic
cylinder function. It is interesting to note that the integrand then
simplifies into the arcsine function and the complete elliptic integral of the
first kind, respectively. Also, a representation is given for the product of a
parabolic cylinder function and a modified Bessel function of the first kind
of order $1/4$.

The remainder of this paper is organized as follows. Section $2$ lists the
definitions and the properties that are used throughout the paper. The
required inverse Laplace transforms are derived in Section $3$. Section $4$
specifies the integral representations for $D_{\nu}\left(  x\right)  D_{\mu
}\left(  y\right)  $ and obtains two expressions for $D_{\nu}\left(
-x\right)  D_{\mu}\left(  y\right)  $. Section $5$ presents some
specializations of the main results.

\section{Definitions and useful properties}

Extensive detail on the parabolic cylinder function, the Laplace transform and
the hypergeometric function can be found in, for instance, \cite{emoth253},
\cite{db07} and \cite{as72}, respectively, such that this section merely lists
the definitions and properties that will be used in the paper.

The parabolic cylinder function of order $\nu$ and argument $z$ is denoted by
$D_{\nu}\left(  z\right)  $. Equation (4) on p. 117 in \cite{emoth253} defines
the parabolic cylinder function in terms of Kummer's confluent hypergeometric
function $\Phi\left(  a;b;z\right)  $ as follows%
\begin{align}
& D_{\nu}\left(  z\right)  =\tag{2.1}\label{dkummer}\\
& 2^{\nu/2}\exp\left(  -\tfrac{1}{4}z^{2}\right)  \left[  \dfrac{\Gamma\left(
1/2\right)  }{\Gamma\left(  \left(  1-\nu\right)  /2\right)  }\Phi\left(
-\dfrac{\nu}{2};\dfrac{1}{2};\dfrac{1}{2}z^{2}\right)  +\dfrac{z}{2^{^{1/2}}%
}\dfrac{\Gamma\left(  -1/2\right)  }{\Gamma\left(  -\nu/2\right)  }\Phi\left(
\dfrac{1-\nu}{2};\dfrac{3}{2};\dfrac{1}{2}z^{2}\right)  \right]  ,\nonumber
\end{align}
\noindent where $\Gamma\left(  \nu\right)  $ denotes the gamma function.
Adding the corresponding relation for $D_{\nu}\left(  -z\right)  $ to
(\ref{dkummer}) gives%
\begin{align}
D_{\nu}\left(  -z\right)  +D_{\nu}\left(  z\right)    & =\dfrac{2^{\left(
\nu+2\right)  /2}\sqrt{\pi}}{\Gamma\left(  \left(  1-\nu\right)  /2\right)
}\exp\left(  -\tfrac{1}{4}z^{2}\right)  \Phi\left(  -\dfrac{\nu}{2};\dfrac
{1}{2};\dfrac{1}{2}z^{2}\right)  ,\tag{2.2a}\label{sum1}\\
D_{\nu}\left(  -z\right)  -D_{\nu}\left(  z\right)    & =\dfrac{z2^{\left(
\nu+3\right)  /2}\sqrt{\pi}}{\Gamma\left(  -\nu/2\right)  }\exp\left(
-\tfrac{1}{4}z^{2}\right)  \Phi\left(  \dfrac{1-\nu}{2};\dfrac{3}{2};\dfrac
{1}{2}z^{2}\right)  ,\tag{2.2b}\label{sum2}%
\end{align}
\noindent see (46:5:3) and (46:5:4) in \cite{oms09}.

The relations in (\ref{sum1}) and (\ref{sum2}) then imply the following two
expressions for $D_{\nu}\left(  -x\right)  D_{\mu}\left(  y\right)  $, namely%
\begin{align}
& D_{\nu}\left(  -x\right)  D_{\mu}\left(  y\right)  =D_{\nu}\left(  x\right)
D_{\mu}\left(  y\right)  +\dfrac{x2^{\left(  \nu+3\right)  /2}\sqrt{\pi}%
\exp\left(  -\tfrac{1}{4}x^{2}\right)  }{\Gamma\left(  -\nu/2\right)  }D_{\mu
}\left(  y\right)  \Phi\left(  \dfrac{1-\nu}{2};\dfrac{3}{2};\dfrac{1}{2}%
x^{2}\right)  ,\nonumber\\
& \tag{2.3a}\label{product1}\\
& D_{\nu}\left(  -x\right)  D_{\mu}\left(  y\right)  =-D_{\nu}\left(
x\right)  D_{\mu}\left(  y\right)  +\dfrac{2^{\left(  \nu+2\right)  /2}%
\sqrt{\pi}\exp\left(  -\tfrac{1}{4}x^{2}\right)  }{\Gamma\left(  \left(
1-\nu\right)  /2\right)  }D_{\mu}\left(  y\right)  \Phi\left(  -\dfrac{\nu}%
{2};\dfrac{1}{2};\dfrac{1}{2}x^{2}\right)  .\nonumber\\
& \tag{2.3b}\label{product2}%
\end{align}
\noindent The following two identities for the parabolic cylinder function
will be used within the below specializations of results%
\begin{align}
& D_{-1/2}\left(  z\right)  =\sqrt{z/(2\pi)}K_{1/4}\left(  z^{2}/4\right)
,\tag{2.4a}\label{d-1/2}\\
& D_{-1}\left(  z\right)  =\exp\left(  z^{2}/4\right)  \sqrt{\pi/2}\text{
erfc}(z/\sqrt{2}),\tag{2.4b}\label{d-1}%
\end{align}
\noindent see p. 326 in \cite{mos66} and (9.254.1) in \cite{gr14}, where
$K_{\nu}\left(  z\right)  $ denotes the modified Bessel function of the second
kind with order $\nu$ and erfc$\left(  z\right)  $ is the complementary error
function. We also use the following relation between $\Phi\left(
a;b;z\right)  $ and $I_{\nu}\left(  z\right)  $, i.e. the modified Bessel
function of the first kind%
\begin{equation}
\Phi\left(  a;2a;z\right)  =\Gamma\left(  a+1/2\right)  \left(  z/4\right)
^{1/2-a}\exp\left(  z/2\right)  I_{a-1/2}\left(  z/2\right)  ,\tag{2.5}%
\label{besseli}%
\end{equation}
\noindent see (5) on p. 579 in \cite{pbm390}.

The convolution theorem will be used to obtain inverse Laplace transforms for
products of two parabolic cylinder functions. As the inverse Laplace
transforms that will be used for $\Phi\left(  a;b;z\right)  $ are defined over
a finite interval rather than for $0\leqslant t<\infty$, the convolution
theorem becomes somewhat more involved. In particular, \cite{pbm492} defines
the Laplace transform of the original functions $f_{1}\left(  t\right)  $ and
$f_{2}\left(  t\right)  $ on p. 578 as%
\begin{align*}
\overline{f}_{1}\left(  p\right)    & =\int_{\alpha_{1}}^{\beta_{1}}%
\exp\left(  -pt\right)  f_{1}\left(  t\right)  dt,\hspace{0.5cm}\beta
_{1}>\alpha_{1},\\
\overline{f}_{2}\left(  p\right)    & =\int_{\alpha_{2}}^{\beta_{2}}%
\exp\left(  -pt\right)  f_{2}\left(  t\right)  dt,\hspace{0.5cm}\beta
_{2}>\alpha_{2},
\end{align*}
\noindent where $\operatorname{Re}p>0$. The convolution theorem then is
specified in \cite{pbm492} as%
\begin{equation}
\overline{f}_{1}\left(  p\right)  \overline{f}_{2}\left(  p\right)
=\int_{\alpha_{1}+\alpha_{2}}^{\beta_{1}+\beta_{2}}\exp\left(  -pt\right)
f_{1}\left(  t\right)  \ast f_{2}\left(  t\right)  dt,\tag{2.6}\label{ilt}%
\end{equation}
\noindent where $f_{1}\left(  t\right)  \ast f_{2}\left(  t\right)  $ is the
convolution of $f_{1}\left(  t\right)  $ and $f_{2}\left(  t\right)  $ that is
obtained from%
\begin{equation}
f_{1}\left(  t\right)  \ast f_{2}\left(  t\right)  =\int_{\max\left(
\alpha_{1};\text{ }t-\beta_{2}\right)  }^{\min\left(  \beta_{1};\text{
}t-\alpha_{2}\right)  }f_{1}\left(  \tau\right)  f_{2}\left(  t-\tau\right)
d\tau.\tag{2.7}\label{convol}%
\end{equation}
\noindent The inverse Laplace transforms for single parabolic cylinder
functions and Kummer's confluent hypergeometric function that will be used
within the convolution theorem (\ref{convol}) are chosen in terms of their
ability to ultimately produce integrands that contain the Gaussian
hypergeometric function $F\left(  a,b;c;z\right)  $. In particular, the latter
function will be obtained in two steps. First, the convolution integral will
be transformed into the integral representation for the Appell function
$F_{1}\left(  a,b_{1},b_{2};c;x,y\right)  $%
\begin{equation}
\dfrac{\Gamma\left(  a\right)  \Gamma\left(  c-a\right)  }{\Gamma\left(
c\right)  }F_{1}\left(  a,b_{1},b_{2};c;x,y\right)  =\int_{0}^{1}%
u^{a-1}\left(  1-u\right)  ^{c-a-1}\left(  1-xu\right)  ^{-b_{1}}\left(
1-yu\right)  ^{-b_{2}}du,\tag{2.8}\label{appell}%
\end{equation}
\noindent for $\operatorname{Re}c>\operatorname{Re}a>0$, see (8.2.5) in
\cite{s66}. Second, the Appell function $F_{1}\left(  a,b_{1},b_{2}%
;c;x,y\right)  $ simplifies into the Gaussian hypergeometric function for
$b_{1}+b_{2}=c$ via%
\begin{equation}
F_{1}\left(  a,b_{1},b_{2};b_{1}+b_{2};x,y\right)  =\left(  1-y\right)
^{-a}F\left(  a,b_{1};b_{1}+b_{2};\dfrac{x-y}{1-y}\right)  ,\tag{2.9}%
\label{gauss}%
\end{equation}
\noindent see (8.3.1.2) in \cite{s66}.

The paper also intensively uses the following (linear transformation) formulas
for $F\left(  a,b;c;z\right)  $%
\begin{align}
F\left(  a,b;c;z\right)    & =F\left(  b,a;c;z\right)  ,\tag{2.10a}%
\label{linear1}\\
& =\left(  1-z\right)  ^{c-a-b}F\left(  c-a,c-b;c;z\right)  ,\tag{2.10b}%
\label{linear2}\\
& =\left(  1-z\right)  ^{-a}F\left(  a,c-b;c;\dfrac{z}{z-1}\right)
,\tag{2.10c}\label{linear3}\\
& =\left(  1-z\right)  ^{-b}F\left(  b,c-a;c;\dfrac{z}{z-1}\right)
,\tag{2.10d}\label{linear4}%
\end{align}
\noindent and the quadratic transformation formula%
\begin{equation}
F\left(  a,b;a+b+1/2;z\right)  =F\left(  2a,2b;a+b+1/2;\left(  1/2\right)
-(1/2)\sqrt{1-z}\right)  ,\tag{2.11}\label{quadratic}%
\end{equation}
\noindent see (15.1.1), (15.3.3)--(15.3.5) and (15.3.22) in \cite{as72}.

Also, the following special cases of $F\left(  a,b;c;z\right)  $ are used%
\begin{align}
& F\left(  a,b;a+b+1/2;x\right)  =2^{a+b-1/2}\Gamma\left(  a+b+1/2\right)
x^{\left(  1-2a-2b\right)  /4}P_{a-b-1/2}^{1/2-a-b}\left(  \sqrt{1-x}\right)
\nonumber\\
& \hspace{1.5cm}\text{for }0<x<1,\tag{2.12a}\label{legendre0}\\
& F\left(  a,b;a-b+1;x\right)  =\Gamma\left(  a-b+1\right)  \left(
1-x\right)  ^{-b}\left(  -x\right)  ^{(b-a)/2}P_{-b}^{b-a}\left(  \left(
1+x\right)  /(1-x)\right)  \nonumber\\
& \hspace{1.5cm}\text{for }-\infty<x<0,\tag{2.12b}\label{legendre1}\\
& F\left(  a,b;b+1;z\right)  =bz^{-b}B\left(  z;b,1-a\right)  ,\tag{2.12c}%
\label{incomplete}\\
& F\left(  a,b;\tfrac{3}{2};x\right)  =-\pi^{-1/2}2^{a+b-7/2}\Gamma\left(
a-1/2\right)  \Gamma\left(  b-1/2\right)  x^{-1/2}\left(  1-x\right)
^{\left(  3-2a-2b\right)  /4}\nonumber\\
& \hspace{1.5cm}\times\left(  P_{a-b-1/2}^{3/2-a-b}\left(  x^{1/2}\right)
-P_{a-b-1/2}^{3/2-a-b}\left(  -x^{1/2}\right)  \right)  \hspace{1.5cm}%
\text{for }0<x<1,\tag{2.12d}\label{legendre2}\\
& F\left(  1/2,1/2;1;z\right)  =\left(  2/\pi\right)  \mathbf{K}\left(
z\right)  ,\tag{2.12e}\label{elliptic}\\
& F\left(  1/2,1/2;3/2;z^{2}\right)  =z^{-1}\arcsin\left(  z\right)
,\tag{2.12f}\label{arcsin}%
\end{align}
\noindent see (15.4.13) and (15.4.15) in \cite{as72}, (28) on p. 455 in
\cite{pbm390} and (15.4.26), (17.3.9) and (15.1.6) in \cite{as72},
respectively. The associated Legendre function of the first kind is denoted by
$P_{b}^{a}\left(  x\right)  $, $B\left(  z;a,b\right)  $ is the incomplete
beta function and $\mathbf{K}\left(  z\right)  $ is the complete elliptic
integral of the first kind.

\section{Three inverse Laplace transforms}

The integral representations for $D_{\nu}\left(  x\right)  D_{\mu}\left(
y\right)  $ and $D_{\nu}\left(  -x\right)  D_{\mu}\left(  y\right)  $ in the
products (\ref{product1}) and (\ref{product2}) will be obtained from three
inverse Laplace transforms.

\subsection{An inverse Laplace transform that contains two parabolic cylinder
functions}

We start from the two inverse Laplace transforms for single parabolic cylinder
functions that are specified in (5) and (6) on p. 290 in \cite{emoti154}%
\begin{align*}
& \Gamma\left(  \nu\right)  \exp\left(  \tfrac{1}{2}ap\right)  D_{-2\nu
}\left(  2^{1/2}a^{1/2}p^{1/2}\right)  =\int_{0}^{\infty}\exp\left(
-pt\right)  2^{-\nu}a^{1/2}t^{\nu-1}\left(  t+a\right)  ^{-\nu-1/2}dt\\
& \hspace{1.5cm}\left[  \operatorname{Re}p>0,\operatorname{Re}\nu>0,\left\vert
\arg a\right\vert <\pi\right]  ,
\end{align*}
\noindent and%
\begin{align*}
& \Gamma\left(  \nu\right)  p^{-1/2}\exp\left(  \tfrac{1}{2}ap\right)
D_{1-2\nu}\left(  2^{1/2}a^{1/2}p^{1/2}\right)  =\int_{0}^{\infty}\exp\left(
-pt\right)  2^{1/2-\nu}t^{\nu-1}\left(  t+a\right)  ^{1/2-\nu}dt\\
& \hspace{1.5cm}\left[  \operatorname{Re}p>0,\operatorname{Re}\nu>0,\left\vert
\arg a\right\vert <\pi\right]  .
\end{align*}
\noindent These inverse transforms, in our notation, can be written as%
\begin{align}
& \Gamma\left(  -\nu/2\right)  \exp\left(  \tfrac{1}{2}xp\right)  D_{\nu
}\left(  2^{1/2}x^{1/2}p^{1/2}\right)  =\int_{0}^{\infty}\exp\left(
-pt\right)  2^{\nu/2}x^{1/2}t^{-\nu/2-1}\left(  t+x\right)  ^{\left(
\nu-1\right)  /2}dt\nonumber\\
& \hspace{1.5cm}\left[  \operatorname{Re}p>0,\operatorname{Re}\nu<0,\left\vert
\arg x\right\vert <\pi\right]  ,\tag{3.1}\label{iltpcf1}%
\end{align}
\noindent and%
\begin{align}
& \Gamma\left(  \left(  1-\mu\right)  /2\right)  p^{-1/2}\exp\left(  \tfrac
{1}{2}yp\right)  D_{\mu}\left(  2^{1/2}y^{1/2}p^{1/2}\right)  =\int
_{0}^{\infty}\exp\left(  -pt\right)  2^{\mu/2}t^{-\left(  \mu+1\right)
/2}\left(  t+y\right)  ^{\mu/2}dt\nonumber\\
& \hspace{1.5cm}\left[  \operatorname{Re}p>0,\operatorname{Re}\mu<1,\left\vert
\arg y\right\vert <\pi\right]  .\tag{3.2}\label{iltpcf2}%
\end{align}
\noindent The original functions $f_{1}\left(  t\right)  $ and $f_{2}\left(
t\right)  $ in the convolution theorem (\ref{ilt}) and (\ref{convol}) are
taken from the inverse transforms (\ref{iltpcf1}) and (\ref{iltpcf2}),
respectively%
\[
f_{1}\left(  t\right)  =2^{\nu/2}x^{1/2}t^{-\nu/2-1}\left(  t+x\right)
^{\left(  \nu-1\right)  /2}\text{ and }f_{2}\left(  t\right)  =2^{\mu
/2}t^{-\left(  \mu+1\right)  /2}\left(  t+y\right)  ^{\mu/2}.
\]
\noindent Given $\beta_{1}=\beta_{2}=\infty$ and $\alpha_{1}=\alpha_{2}=0$ in
(\ref{ilt}), $\max\left(  \alpha_{1};t-\beta_{2}\right)  =0$ and $\min\left(
\beta_{1};t-\alpha_{2}\right)  =t$ such that the convolution integral is%
\[
f_{1}\left(  t\right)  \ast f_{2}\left(  t\right)  =\int_{0}^{t}2^{\nu
/2}x^{1/2}\tau^{-\nu/2-1}\left(  \tau+x\right)  ^{\left(  \nu-1\right)
/2}2^{\mu/2}\left(  t-\tau\right)  ^{-\left(  \mu+1\right)  /2}\left(  \left(
t-\tau\right)  +y\right)  ^{\mu/2}d\tau.
\]
\noindent The substitution $\tau=tu$ gives%
\begin{align*}
& f_{1}\left(  t\right)  \ast f_{2}\left(  t\right)  =2^{\left(  \nu
+\mu\right)  /2}x^{\nu/2}t^{-\left(  1+\nu+\mu\right)  /2}\left(  y+t\right)
^{\mu/2}\\
& \hspace{1cm}\times\int_{0}^{1}u^{-\nu/2-1}\left(  1-u\right)  ^{-\left(
1+\mu\right)  /2}\left(  1+\dfrac{t}{x}u\right)  ^{\left(  \nu-1\right)
/2}\left(  1-\dfrac{t}{t+y}u\right)  ^{\mu/2}du.
\end{align*}
\noindent The definition (\ref{appell}) allows to rewrite this result as%
\begin{align*}
& f_{1}\left(  t\right)  \ast f_{2}\left(  t\right)  =2^{\left(  \nu
+\mu\right)  /2}x^{\nu/2}t^{-\left(  1+\nu+\mu\right)  /2}\left(  y+t\right)
^{\mu/2}\dfrac{\Gamma\left(  -\nu/2\right)  \Gamma\left(  \left(
1-\mu\right)  /2\right)  }{\Gamma\left(  \left(  1-\nu-\mu\right)  /2\right)
}\\
& \hspace{1cm}\times F_{1}\left(  -\tfrac{\nu}{2},\tfrac{1-\nu}{2},-\tfrac
{\mu}{2};\tfrac{1}{2}\left(  1-\nu-\mu\right)  ;-\dfrac{t}{x},\dfrac{t}%
{t+y}\right)  .
\end{align*}
\noindent As $b_{1}+b_{2}=c$ in (\ref{appell}), we have%
\begin{align*}
& f_{1}\left(  t\right)  \ast f_{2}\left(  t\right)  =2^{\left(  \nu
+\mu\right)  /2}\left(  xy\right)  ^{\nu/2}t^{-\left(  1+\nu+\mu\right)
/2}\left(  y+t\right)  ^{\left(  \mu-\nu\right)  /2}\dfrac{\Gamma\left(
-\nu/2\right)  \Gamma\left(  \left(  1-\mu\right)  /2\right)  }{\Gamma\left(
\left(  1-\nu-\mu\right)  /2\right)  }\\
& \hspace{1cm}\times F\left(  -\tfrac{\nu}{2},\tfrac{1-\nu}{2};\tfrac{1}%
{2}\left(  1-\nu-\mu\right)  ;-\dfrac{t\left(  x+y+t\right)  }{xy}\right)  ,
\end{align*}
\noindent which via (\ref{linear3}) is rewritten into%
\begin{align*}
& f_{1}\left(  t\right)  \ast f_{2}\left(  t\right)  =2^{\left(  \nu
+\mu\right)  /2}t^{-\left(  1+\nu+\mu\right)  /2}\left(  y+t\right)  ^{\mu
/2}\left(  x+t\right)  ^{\nu/2}\dfrac{\Gamma\left(  -\nu/2\right)
\Gamma\left(  \left(  1-\mu\right)  /2\right)  }{\Gamma\left(  \left(
1-\nu-\mu\right)  /2\right)  }\\
& \hspace{1cm}\times F\left(  -\tfrac{\nu}{2},-\tfrac{\mu}{2};\tfrac{1}%
{2}\left(  1-\nu-\mu\right)  ;\dfrac{t\left(  x+y+t\right)  }{\left(
x+t\right)  \left(  y+t\right)  }\right)  .
\end{align*}
\noindent Using (\ref{ilt}) with $\alpha_{1}+\alpha_{2}=0$ and $\beta
_{1}+\beta_{2}=\infty$ then gives the following inverse Laplace transform%
\begin{align}
& p^{-1/2}\exp\left(  \tfrac{1}{2}p\left(  x+y\right)  \right)  D_{\nu}\left(
2^{1/2}x^{1/2}p^{1/2}\right)  D_{\mu}\left(  2^{1/2}y^{1/2}p^{1/2}\right)
=\nonumber\\
& \dfrac{2^{\left(  \nu+\mu\right)  /2}}{\Gamma\left(  \left(  1-\nu
-\mu\right)  /2\right)  }\int_{0}^{\infty}\exp\left(  -pt\right)  t^{-\left(
1+\nu+\mu\right)  /2}\left(  y+t\right)  ^{\mu/2}\left(  x+t\right)  ^{\nu
/2}\nonumber\\
& \hspace{1cm}\times F\left(  -\tfrac{\nu}{2},-\tfrac{\mu}{2};\tfrac{1}%
{2}\left(  1-\nu-\mu\right)  ;\dfrac{t\left(  x+y+t\right)  }{\left(
x+t\right)  \left(  y+t\right)  }\right)  dt\nonumber\\
& \hspace{1.5cm}\left[  \operatorname{Re}p>0,\operatorname{Re}\nu
<0,\operatorname{Re}\mu<1,\left\vert \arg x\right\vert <\pi,\left\vert \arg
y\right\vert <\pi\right]  .\tag{3.3}\label{iltpcfpcf}%
\end{align}

\subsection{Two inverse Laplace transforms that contain a parabolic cylinder
function and a Kummer confluent hypergeometric function}

Evaluating (\ref{product1}) requires an inverse Laplace transform that
contains $D_{\mu}\left(  y\right)  \Phi\left(  \left(  1-\nu\right)
/2;3/2;x^{2}/2\right)  $. The inverse Laplace transform for $D_{\mu}\left(
y\right)  $ is taken from (\ref{iltpcf1}) in which $\nu$ and $x$ are
interchanged for $\mu$ and $y$, respectively. The inverse transform for
$\Phi\left(  \left(  1-\nu\right)  /2;3/2;x^{2}/2\right)  $ is obtained from
the inverse Laplace transform (3.33.2.2) in \cite{pbm592}, namely%
\begin{align}
& \exp\left(  -xp\right)  \Phi\left(  a;b;xp\right)  =\dfrac{x^{1-b}%
\Gamma\left(  b\right)  }{\Gamma\left(  b-a\right)  \Gamma\left(  a\right)
}\int_{0}^{x}\exp\left(  -pt\right)  t^{b-a-1}\left(  x-t\right)
^{a-1}dt\nonumber\\
& \hspace{1.5cm}\left[  \operatorname{Re}p>0,\operatorname{Re}%
b>\operatorname{Re}a>0,x>0\right]  .\tag{3.4}\label{iltkum}%
\end{align}
\noindent Specializing the latter for $a=\left(  1-\nu\right)  /2$ and $b=3/2$
with rescaling of $x$ to $x^{2}/2$ gives%
\begin{align}
& x\sqrt{2/\pi}\Gamma\left(  1+\nu/2\right)  \Gamma\left(  \left(
1-\nu\right)  /2\right)  \exp\left(  -\tfrac{1}{2}x^{2}p\right)  \Phi\left(
\tfrac{1-\nu}{2};\tfrac{3}{2};\tfrac{1}{2}px^{2}\right)  =\nonumber\\
& \hspace{1cm}\int_{0}^{\tfrac{1}{2}x^{2}}\exp\left(  -pt\right)  t^{\nu
/2}\left(  \tfrac{1}{2}x^{2}-t\right)  ^{-\left(  1+\nu\right)  /2}%
dt\nonumber\\
& \hspace{1.5cm}\left[  \operatorname{Re}p>0,-2<\operatorname{Re}%
\nu<1,x>0\right]  .\tag{3.5}\label{iltkum1}%
\end{align}
\noindent The original functions $f_{1}\left(  t\right)  $ and $f_{2}\left(
t\right)  $ are taken from the inverse Laplace transforms of the
aforementioned parabolic cylinder function and (\ref{iltkum1}), respectively%
\[
f_{1}\left(  t\right)  =2^{\mu/2}y^{1/2}t^{-\mu/2-1}\left(  t+y\right)
^{\left(  \mu-1\right)  /2}\text{ and }f_{2}\left(  t\right)  =t^{\nu
/2}\left(  \tfrac{1}{2}x^{2}-t\right)  ^{-\left(  1+\nu\right)  /2}.
\]
\noindent The integration limits in (\ref{ilt}) and (\ref{convol}) are
$\beta_{1}=\infty,\beta_{2}=x^{2}/2,\alpha_{1}=\alpha_{2}=0$ and the
convolution integral is%
\begin{equation}%
\begin{tabular}
[c]{lll}%
$f_{1}\left(  t\right)  \ast f_{2}\left(  t\right)  $ & $=%
{\displaystyle\int\nolimits_{0}^{t}}
f_{1}\left(  \tau\right)  f_{2}\left(  t-\tau\right)  d\tau$ & $\hspace
{1cm}t<\tfrac{1}{2}x^{2},$\\
&  & \\
& $=%
{\displaystyle\int\nolimits_{t-\tfrac{1}{2}x^{2}}^{t}}
f_{1}\left(  \tau\right)  f_{2}\left(  t-\tau\right)  d\tau$ & $\hspace
{1cm}t>\tfrac{1}{2}x^{2}.$%
\end{tabular}
\tag{3.6}\label{convol1}%
\end{equation}
\noindent The convolution integral for $t<x^{2}/2$ then is%
\[
f_{1}\left(  t\right)  \ast f_{2}\left(  t\right)  =\int_{0}^{t}2^{\mu
/2}y^{1/2}\tau^{-\mu/2-1}\left(  \tau+y\right)  ^{\left(  \mu-1\right)
/2}\left(  t-\tau\right)  ^{\nu/2}\left(  \tfrac{1}{2}x^{2}-\left(
t-\tau\right)  \right)  ^{-\left(  1+\nu\right)  /2}d\tau.
\]
\noindent Using the substitution $\tau=tu$ and employing steps akin to those
used within the previous inverse Laplace transform yield%
\begin{align*}
& f_{1}\left(  t\right)  \ast f_{2}\left(  t\right)  =2^{\left(  1+\nu
+\mu\right)  /2}t^{\left(  \nu-\mu\right)  /2}\left(  x^{2}-2t\right)
^{-\left(  1+\nu\right)  /2}\left(  y+t\right)  ^{\mu/2}\dfrac{\Gamma\left(
-\mu/2\right)  \Gamma\left(  1+\nu/2\right)  }{\Gamma\left(  1+\left(  \nu
-\mu\right)  /2\right)  }\\
& \hspace{1cm}\times F\left(  -\tfrac{\mu}{2},\tfrac{\nu+1}{2};1+\tfrac{1}%
{2}\left(  \nu-\mu\right)  ;\dfrac{t\left(  x^{2}-2y-2t\right)  }{\left(
x^{2}-2t\right)  \left(  y+t\right)  }\right)  \\
& \hspace{1.5cm}\left[  t<\tfrac{1}{2}x^{2}\right]  .
\end{align*}
\noindent The convolution integral for $t>x^{2}/2$ is given by%
\[
f_{1}\left(  t\right)  \ast f_{2}\left(  t\right)  =\int_{t-\tfrac{1}{2}x^{2}%
}^{t}2^{\mu/2}y^{1/2}\tau^{-\mu/2-1}\left(  \tau+y\right)  ^{\left(
\mu-1\right)  /2}\left(  t-\tau\right)  ^{\nu/2}\left(  \tfrac{1}{2}%
x^{2}-\left(  t-\tau\right)  \right)  ^{-\left(  1+\nu\right)  /2}d\tau.
\]
\noindent The substitution $\tau=s-x^{2}/2+t$ yields%
\begin{align*}
& f_{1}\left(  t\right)  \ast f_{2}\left(  t\right)  =\\
& \hspace{1cm}2^{\mu/2}y^{1/2}\int_{0}^{\tfrac{1}{2}x^{2}}\left(  t-\tfrac
{1}{2}x^{2}+s\right)  ^{-\mu/2-1}\left(  t+y-\tfrac{1}{2}x^{2}+s\right)
^{\left(  \mu-1\right)  /2}\left(  \tfrac{1}{2}x^{2}-s\right)  ^{\nu
/2}s^{-\left(  1+\nu\right)  /2}ds.
\end{align*}
\noindent Using the substitution $s=x^{2}u/2$ and employing the above steps
toward $F\left(  a,b;c;z\right)  $ with the use of the linear transformation
formula (\ref{linear3}) give%
\begin{align*}
& f_{1}\left(  t\right)  \ast f_{2}\left(  t\right)  =\dfrac{2^{\left(
3+\nu+\mu\right)  /2}}{\sqrt{\pi}}xy^{1/2}t^{\left(  \nu-1\right)  /2}\left(
2t-x^{2}\right)  ^{-\left(  1+\nu+\mu\right)  /2}\left(  2y-x^{2}+2t\right)
^{\left(  \mu-1\right)  /2}\\
& \hspace{1cm}\times\Gamma\left(  \left(  1-\nu\right)  /2\right)
\Gamma\left(  1+\nu/2\right)  F\left(  \tfrac{1-\nu}{2},\tfrac{1-\mu}%
{2};\tfrac{3}{2};\dfrac{x^{2}y}{t\left(  2y-x^{2}+2t\right)  }\right)  \\
& \hspace{1.5cm}\left[  t>\tfrac{1}{2}x^{2}\right]  .
\end{align*}
\noindent The inverse Laplace transform (\ref{ilt}) with $\alpha_{1}%
+\alpha_{2}=0$ and $\beta_{1}+\beta_{2}=\infty$ then is%
\begin{align}
& \exp\left(  \tfrac{1}{2}p\left(  y-x^{2}\right)  \right)  D_{\mu}\left(
2^{1/2}y^{1/2}p^{1/2}\right)  \Phi\left(  \tfrac{1-\nu}{2};\tfrac{3}{2}%
;\tfrac{1}{2}px^{2}\right)  =\nonumber\\
& \hspace{1cm}\dfrac{2^{\left(  \nu+\mu\right)  /2}\sqrt{\pi}x^{-1}}%
{\Gamma\left(  \left(  1-\nu\right)  /2\right)  \Gamma\left(  1+\left(
\nu-\mu\right)  /2\right)  }\int_{0}^{\tfrac{1}{2}x^{2}}\exp\left(
-pt\right)  t^{\left(  \nu-\mu\right)  /2}\left(  x^{2}-2t\right)  ^{-\left(
\nu+1\right)  /2}\left(  y+t\right)  ^{\mu/2}\nonumber\\
& \hspace{1cm}\times F\left(  -\tfrac{\mu}{2},\tfrac{\nu+1}{2};1+\tfrac
{\nu-\mu}{2};\dfrac{t\left(  x^{2}-2y-2t\right)  }{\left(  x^{2}-2t\right)
\left(  y+t\right)  }\right)  dt\nonumber\\
& +\dfrac{2^{\left(  2+\nu+\mu\right)  /2}y^{1/2}}{\Gamma\left(
-\mu/2\right)  }\int_{\tfrac{1}{2}x^{2}}^{\infty}\exp\left(  -pt\right)
t^{\left(  \nu-1\right)  /2}\left(  2t-x^{2}\right)  ^{-\left(  1+\nu
+\mu\right)  /2}\left(  2y-x^{2}+2t\right)  ^{\left(  \mu-1\right)
/2}\nonumber\\
& \hspace{1cm}\times F\left(  \tfrac{1-\nu}{2},\tfrac{1-\mu}{2};\tfrac{3}%
{2};\dfrac{x^{2}y}{t\left(  2y-x^{2}+2t\right)  }\right)  dt\nonumber\\
& \hspace{1.5cm}\left[  \operatorname{Re}p>0,-2<\operatorname{Re}%
\nu<1,\operatorname{Re}\mu<0,x>0,\left\vert \arg y\right\vert <\pi
,y\neq0\right]  .\tag{3.7}\label{iltpcfkum1}%
\end{align}
\noindent The inverse Laplace transform for $D_{\mu}\left(  y\right)
\Phi\left(  -\nu/2;1/2;x^{2}/2\right)  $ that is later to be used within the
product (\ref{product2}) starts from the inverse transform (\ref{iltpcf2}) and
the specialization of (\ref{iltkum}) for $a=-\nu/2$ and $b=1/2$%
\begin{align}
& x^{-1}\sqrt{2/\pi}\Gamma\left(  \left(  1+\nu\right)  /2\right)
\Gamma\left(  -\nu/2\right)  \exp\left(  -\tfrac{1}{2}x^{2}p\right)
\Phi\left(  -\tfrac{\nu}{2};\tfrac{1}{2};\tfrac{1}{2}px^{2}\right)
=\nonumber\\
& \hspace{1cm}\int_{0}^{\tfrac{1}{2}x^{2}}\exp\left(  -pt\right)  t^{\left(
\nu-1\right)  /2}\left(  \tfrac{1}{2}x^{2}-t\right)  ^{-\nu/2-1}dt\nonumber\\
& \hspace{1.5cm}\left[  \operatorname{Re}p>0,-1<\operatorname{Re}%
\nu<0,x>0\right]  .\tag{3.8}\label{iltkum2}%
\end{align}
\noindent Taking $f_{1}\left(  t\right)  $ from the parabolic cylinder
function and $f_{2}\left(  t\right)  $ from the Kummer confluent
hypergeometric function implies%
\[
f_{1}\left(  t\right)  =2^{\mu/2}t^{-\left(  1+\mu\right)  /2}\left(
t+y\right)  ^{\mu/2}\text{ and }f_{2}\left(  t\right)  =t^{\left(
\nu-1\right)  /2}\left(  \tfrac{1}{2}x^{2}-t\right)  ^{-\nu/2-1}.
\]
\noindent The convolution integral is identical to (\ref{convol1}) as the
integration limits again are $\beta_{1}=\infty,\beta_{2}=x^{2}/2$ and
$\alpha_{1}=\alpha_{2}=0$. Repeating the above steps then gives the following
inverse Laplace transform%
\begin{align}
& p^{-1/2}\exp\left(  \tfrac{1}{2}p\left(  y-x^{2}\right)  \right)  D_{\mu
}\left(  2^{1/2}y^{1/2}p^{1/2}\right)  \Phi\left(  -\tfrac{\nu}{2};\tfrac
{1}{2};\tfrac{1}{2}px^{2}\right)  =\nonumber\\
& \dfrac{xy^{1/2}\sqrt{\pi}2^{\left(  1+\nu+\mu\right)  /2}}{\Gamma\left(
-\nu/2\right)  \Gamma\left(  1+\left(  \nu-\mu\right)  /2\right)  }\int
_{0}^{\tfrac{1}{2}x^{2}}\exp\left(  -pt\right)  t^{\left(  \nu-\mu\right)
/2}\left(  x^{2}-2t\right)  ^{-\nu/2-1}\left(  y+t\right)  ^{\left(
\mu-1\right)  /2}\nonumber\\
& \hspace{1cm}\times F\left(  \tfrac{\left(  1-\mu\right)  }{2},1+\tfrac{\nu
}{2};1+\tfrac{\nu-\mu}{2};\dfrac{t\left(  x^{2}-2y-2t\right)  }{\left(
x^{2}-2t\right)  \left(  y+t\right)  }\right)  dt\nonumber\\
& +\dfrac{2^{\left(  1+\nu+\mu\right)  /2}}{\Gamma\left(  \left(
1-\mu\right)  /2\right)  }\int_{\tfrac{1}{2}x^{2}}^{\infty}\exp\left(
-pt\right)  t^{\nu/2}\left(  2t-x^{2}\right)  ^{-\left(  1+\nu+\mu\right)
/2}\left(  2y-x^{2}+2t\right)  ^{\mu/2}\nonumber\\
& \hspace{1cm}\times F\left(  -\tfrac{\nu}{2},-\tfrac{\mu}{2};\tfrac{1}%
{2};\dfrac{x^{2}y}{t\left(  2y-x^{2}+2t\right)  }\right)  dt\nonumber\\
& \hspace{1.5cm}\left[  \operatorname{Re}p>0,-1<\operatorname{Re}%
\nu<0,\operatorname{Re}\mu<1,x>0,\left\vert \arg y\right\vert <\pi\right]
.\tag{3.9}\label{iltpcfkum2}%
\end{align}

\section{Integral representations for products of parabolic cylinder
functions}

Setting $p=1$ in the inverse Laplace transform (\ref{iltpcfpcf}) and rescaling
$2^{1/2}x^{1/2}$ and $2^{1/2}y^{1/2}$ to $x$ and $y$, respectively, gives the
following integral representation for $D_{\nu}\left(  x\right)  D_{\mu}\left(
y\right)  $%
\begin{align}
& D_{\nu}\left(  x\right)  D_{\mu}\left(  y\right)  =\dfrac{\exp\left(
-\tfrac{1}{4}\left(  x^{2}+y^{2}\right)  \right)  }{\Gamma\left(  \left(
1-\nu-\mu\right)  /2\right)  }\int_{0}^{\infty}\exp\left(  -t\right)
t^{-\left(  1+\nu+\mu\right)  /2}\left(  x^{2}+2t\right)  ^{\nu/2}\left(
y^{2}+2t\right)  ^{\mu/2}\nonumber\\
& \hspace{1cm}\times F\left(  -\tfrac{\nu}{2},-\tfrac{\mu}{2};\tfrac{1}%
{2}\left(  1-\nu-\mu\right)  ;\dfrac{2t\left(  x^{2}+y^{2}+2t\right)
}{\left(  x^{2}+2t\right)  \left(  y^{2}+2t\right)  }\right)  dt\nonumber\\
& \hspace{1.5cm}\left[  \operatorname{Re}\nu<0,\operatorname{Re}%
\mu<1,x\geqslant0,y\geqslant0\right]  .\tag{4.1}\label{representation1}%
\end{align}
\noindent Note that the linear transformation formulas (\ref{linear2}%
)--(\ref{linear4}) and the quadratic transformation formula (\ref{quadratic})
generate alternative expressions for $D_{\nu}\left(  x\right)  D_{\mu}\left(
y\right)  $ in which the integrands differ but still evolve around the
Gaussian hypergeometric function. The identity (\ref{legendre0}), on the
contrary, allows to rephrase the integrand of (\ref{representation1}) in terms
of the associated Legendre function of the first kind%
\begin{align}
& D_{\nu}\left(  x\right)  D_{\mu}\left(  y\right)  =\exp\left(  -\tfrac{1}%
{4}\left(  x^{2}+y^{2}\right)  \right)  \int_{0}^{\infty}\exp\left(
-t\right)  \left(  2t\right)  ^{-\left(  1+\nu+\mu\right)  /4}\left(
x^{2}+2t\right)  ^{\left(  \nu-\mu-1\right)  /4}\nonumber\\
& \hspace{1cm}\times\left(  y^{2}+2t\right)  ^{\left(  \mu-\nu-1\right)
/4}\left(  x^{2}+y^{2}+2t\right)  ^{\left(  1+\nu+\mu\right)  /4}P_{\left(
\mu-\nu-1\right)  /2}^{\left(  1+\nu+\mu\right)  /2}\left(  \dfrac{xy}%
{\sqrt{\left(  x^{2}+2t\right)  \left(  y^{2}+2t\right)  }}\right)
dt\nonumber\\
& \hspace{1.5cm}\left[  \operatorname{Re}\nu<0,\operatorname{Re}%
\mu<1,x\geqslant0,y\geqslant0\right]  .\tag{4.2}\label{representation1b}%
\end{align}
\noindent The first integral representation for $D_{\nu}\left(  -x\right)
D_{\mu}\left(  y\right)  $ is obtained by setting $p=1$ and rescaling
$2^{1/2}y^{1/2}\rightarrow y$ in the inverse Laplace transform
(\ref{iltpcfkum1}) and plugging the outcome together with
(\ref{representation1}) into the product (\ref{product1})%
\begin{align}
& D_{\nu}\left(  -x\right)  D_{\mu}\left(  y\right)  =\dfrac{\exp\left(
-\tfrac{1}{4}\left(  x^{2}+y^{2}\right)  \right)  }{\Gamma\left(  \left(
1-\nu-\mu\right)  /2\right)  }\int_{0}^{\infty}\exp\left(  -t\right)
t^{-\left(  1+\nu+\mu\right)  /2}\left(  x^{2}+2t\right)  ^{\nu/2}\left(
y^{2}+2t\right)  ^{\mu/2}\nonumber\\
& \hspace{1cm}\times F\left(  -\tfrac{\nu}{2},-\tfrac{\mu}{2};\tfrac{1}%
{2}\left(  1-\nu-\mu\right)  ;\dfrac{2t\left(  x^{2}+y^{2}+2t\right)
}{\left(  x^{2}+2t\right)  \left(  y^{2}+2t\right)  }\right)  dt\nonumber\\
& +\dfrac{\sqrt{2\pi}\exp\left(  \tfrac{1}{4}\left(  x^{2}-y^{2}\right)
\right)  }{\Gamma\left(  -\nu\right)  \Gamma\left(  1+\left(  \nu-\mu\right)
/2\right)  }\int_{0}^{\tfrac{1}{2}x^{2}}\exp\left(  -t\right)  t^{\left(
\nu-\mu\right)  /2}\left(  x^{2}-2t\right)  ^{-\left(  \nu+1\right)
/2}\left(  y^{2}+2t\right)  ^{\mu/2}\nonumber\\
& \hspace{1cm}\times F\left(  -\tfrac{\mu}{2},\tfrac{\nu+1}{2};1+\tfrac
{\nu-\mu}{2};\dfrac{2t\left(  x^{2}-y^{2}-2t\right)  }{\left(  x^{2}%
-2t\right)  \left(  y^{2}+2t\right)  }\right)  dt\nonumber\\
& +\dfrac{xy\sqrt{\pi}2^{2+\nu+\mu/2}\exp\left(  \tfrac{1}{4}\left(
x^{2}-y^{2}\right)  \right)  }{\Gamma\left(  -\nu/2\right)  \Gamma\left(
-\mu/2\right)  }\int_{\tfrac{1}{2}x^{2}}^{\infty}\exp\left(  -t\right)
t^{\left(  \nu-1\right)  /2}\left(  2t-x^{2}\right)  ^{-\left(  1+\nu
+\mu\right)  /2}\nonumber\\
& \hspace{1cm}\times\left(  y^{2}-x^{2}+2t\right)  ^{\left(  \mu-1\right)
/2}F\left(  \tfrac{1-\nu}{2},\tfrac{1-\mu}{2};\tfrac{3}{2};\dfrac{x^{2}y^{2}%
}{2t\left(  y^{2}-x^{2}+2t\right)  }\right)  dt\nonumber\\
& \hspace{1.5cm}\left[  -2<\operatorname{Re}\nu<0,\operatorname{Re}%
\mu<0,x\geqslant0,y>0\right]  .\tag{4.3}\label{representation2}%
\end{align}
\noindent Setting $p=1$ and rescaling $2^{1/2}y^{1/2}\rightarrow y$ in
(\ref{iltpcfkum2}) and subsequently inserting the result together with
(\ref{representation1}) into the product (\ref{product2}) gives a second
integral representation for $D_{\nu}\left(  -x\right)  D_{\mu}\left(
y\right)  $%
\begin{align}
& D_{\nu}\left(  -x\right)  D_{\mu}\left(  y\right)  =-\dfrac{\exp\left(
-\tfrac{1}{4}\left(  x^{2}+y^{2}\right)  \right)  }{\Gamma\left(  \left(
1-\nu-\mu\right)  /2\right)  }\int_{0}^{\infty}\exp\left(  -t\right)
t^{-\left(  1+\nu+\mu\right)  /2}\left(  x^{2}+2t\right)  ^{\nu/2}\left(
y^{2}+2t\right)  ^{\mu/2}\nonumber\\
& \hspace{1cm}\times F\left(  -\tfrac{\nu}{2},-\tfrac{\mu}{2};\tfrac{1}%
{2}\left(  1-\nu-\mu\right)  ;\dfrac{2t\left(  x^{2}+y^{2}+2t\right)
}{\left(  x^{2}+2t\right)  \left(  y^{2}+2t\right)  }\right)  dt\nonumber\\
& +\dfrac{xy\sqrt{2\pi}\exp\left(  \tfrac{1}{4}\left(  x^{2}-y^{2}\right)
\right)  }{\Gamma\left(  -\nu\right)  \Gamma\left(  1+\left(  \nu-\mu\right)
/2\right)  }\int_{0}^{\tfrac{1}{2}x^{2}}\exp\left(  -t\right)  t^{\left(
\nu-\mu\right)  /2}\left(  x^{2}-2t\right)  ^{-\nu/2-1}\left(  y^{2}%
+2t\right)  ^{\left(  \mu-1\right)  /2}\nonumber\\
& \hspace{1cm}\times F\left(  \tfrac{1-\mu}{2},1+\tfrac{\nu}{2};1+\tfrac
{\nu-\mu}{2};\dfrac{2t\left(  x^{2}-y^{2}-2t\right)  }{\left(  x^{2}%
-2t\right)  \left(  y^{2}+2t\right)  }\right)  dt\nonumber\\
& +\dfrac{2^{\left(  3+2\nu+\mu\right)  /2}\sqrt{\pi}\exp\left(  \tfrac{1}%
{4}\left(  x^{2}-y^{2}\right)  \right)  }{\Gamma\left(  \left(  1-\nu\right)
/2\right)  \Gamma\left(  \left(  1-\mu\right)  /2\right)  }\int_{\tfrac{1}%
{2}x^{2}}^{\infty}\exp\left(  -t\right)  t^{\nu/2}\left(  2t-x^{2}\right)
^{-\left(  1+\nu+\mu\right)  /2}\nonumber\\
& \hspace{1cm}\times\left(  y^{2}-x^{2}+2t\right)  ^{\mu/2}F\left(
-\tfrac{\nu}{2},-\tfrac{\mu}{2};\tfrac{1}{2};\dfrac{x^{2}y^{2}}{2t\left(
y^{2}-x^{2}+2t\right)  }\right)  dt\nonumber\\
& \hspace{1.5cm}\left[  -1<\operatorname{Re}\nu<0,\operatorname{Re}%
\mu<1,x>0,y>0\right]  .\tag{4.4}\label{representation3}%
\end{align}

\noindent The latter two expressions for $D_{\nu}\left(  -x\right)  D_{\mu
}\left(  y\right)  $ are complementary as they are valid in different regions
for the real parts of $\nu$ and $\mu$.

\section{Some specializations}

The above results allow for various specializations in which other special
functions emerge both on the left-hand side as well as in the integrand. The
first two examples employ the special values for the parabolic cylinder
function in the identities (\ref{d-1/2}) and (\ref{d-1}). For $\nu=\mu=-1/2$,
(\ref{representation1}) simplifies into an integral representation for the
product of two modified Bessel functions of the second kind of order $1/4$.
The Gaussian hypergeometric function in the integrand then reduces into
$F\left(  1/4,1/4;1;z\right)  =\left(  2/\pi\right)  \mathbf{K}\left(
1/2\left(  1-\sqrt{1-z}\right)  \right)  $ as can be seen from
(\ref{quadratic}) and (\ref{elliptic})%
\begin{align*}
& K_{1/4}\left(  x\right)  K_{1/4}\left(  y\right)  =2^{1/2}x^{-1/4}%
y^{-1/4}\exp\left(  -\left(  x+y\right)  \right)  \int_{0}^{\infty}\exp\left(
-t\right)  \left(  2x+t\right)  ^{-1/4}\left(  2y+t\right)  ^{-1/4}\\
& \hspace{1cm}\times\mathbf{K}\left(  \tfrac{1}{2}-\sqrt{\dfrac{xy}{\left(
2x+t\right)  \left(  2y+t\right)  }}\right)  dt\\
& \hspace{1.5cm}\left[  x>0,y>0\right]  .
\end{align*}
\noindent Setting $\nu=\mu=-1$ in (\ref{representation1}) and using
(\ref{d-1}) yields an integral representation for the product of two
complementary error functions for which the integrand contains the arcsine
function in view of identity (\ref{arcsin})%
\begin{align*}
& \text{erfc}\left(  x\right)  \text{erfc}\left(  y\right)  =2\exp\left(
-\left(  x^{2}+y^{2}\right)  \right)  \pi^{-3/2}\int_{0}^{\infty}\exp\left(
-t\right)  \left(  x^{2}+y^{2}+t\right)  ^{-1/2}\\
& \hspace{1cm}\times\arcsin\left(  \sqrt{\dfrac{t\left(  x^{2}+y^{2}+t\right)
}{\left(  x^{2}+t\right)  \left(  y^{2}+t\right)  }}\right)  dt\\
& \hspace{1.5cm}\left[  x\geqslant0,y\geqslant0\right]  .
\end{align*}
\noindent The next example illustrates an integral representation for the
product of a parabolic cylinder function with a special function that follows
from choosing a special value for $\Phi\left(  a;b;z\right)  $ in the inverse
Laplace transform (\ref{iltpcfkum1}). First, the latter inverse transform
gives the following integral representation%
\begin{align}
& D_{\nu}\left(  x\right)  \Phi\left(  \tfrac{1-\mu}{2};\tfrac{3}{2};y\right)
=\dfrac{2^{-1}\sqrt{\pi}\exp\left(  -\tfrac{1}{4}\left(  x^{2}-4y\right)
\right)  y^{-1/2}}{\Gamma\left(  \left(  1-\mu\right)  /2\right)
\Gamma\left(  1+\left(  \mu-\nu\right)  /2\right)  }\int_{0}^{y}\exp\left(
-t\right)  t^{\left(  \mu-\nu\right)  /2}\left(  x^{2}+2t\right)  ^{\nu
/2}\nonumber\\
& \hspace{1cm}\times\left(  y-t\right)  ^{-\left(  \mu+1\right)  /2}F\left(
-\tfrac{\nu}{2},\tfrac{\mu+1}{2};1+\tfrac{\mu-\nu}{2};\dfrac{t\left(
x^{2}-2y+2t\right)  }{\left(  x^{2}+2t\right)  \left(  t-y\right)  }\right)
dt\nonumber\\
& +\dfrac{x\exp\left(  -\tfrac{1}{4}\left(  x^{2}-4y\right)  \right)  }%
{\Gamma\left(  -\nu/2\right)  }\int_{y}^{\infty}\exp\left(  -t\right)
t^{\left(  \mu-1\right)  /2}\left(  t-y\right)  ^{-\left(  1+\nu+\mu\right)
/2}\left(  x^{2}-2y+2t\right)  ^{\left(  \nu-1\right)  /2}\nonumber\\
& \hspace{1cm}\times F\left(  \tfrac{1-\mu}{2},\tfrac{1-\nu}{2};\tfrac{3}%
{2};\dfrac{x^{2}y}{t\left(  x^{2}-2y+2t\right)  }\right)  dt\nonumber\\
& \hspace{1.5cm}\left[  \operatorname{Re}\nu<0,-2<\operatorname{Re}%
\mu<1,x>0,y>0\right]  .\tag{5.1}\label{representation4}%
\end{align}
\noindent The Kummer confluent hypergeometric function in this representation
is then linked to the modified Bessel function of the first kind with order
$1/4$ for $\mu=-1/2$ in view of the identity (\ref{besseli}). The Gaussian
hypergeometric function in the first integral in (\ref{representation4}) can
be expressed in terms of $P_{b}^{a}\left(  x\right)  $ via relation
(\ref{legendre1}). The second integrand contains, in stylized form, $F\left(
a,b;3/2;z\right)  $ for which the identity (\ref{legendre2}) can be used%
\begin{align*}
& D_{\nu}\left(  x\right)  I_{1/4}(y)=\dfrac{\exp\left(  -\tfrac{1}{4}\left(
x^{2}-4y\right)  \right)  }{x^{1/2}y^{1/2}\sqrt{2\pi}}\int_{0}^{2y}\exp\left(
-t\right)  \left(  \frac{\left(  x^{2}+2t\right)  \left(  x^{2}-4y+2t\right)
}{t\left(  2y-t\right)  }\right)  ^{(1+2\nu)/8}\\
& \hspace{1cm}\times P_{-1/4}^{\left(  1+2\nu\right)  /4}\left(  \frac
{x^{2}y-t\left(  x^{2}-4y+2t\right)  }{x^{2}y}\right)  dt\\
& +\dfrac{2^{-1-\left(  \nu/2\right)  }\exp\left(  -\tfrac{1}{4}\left(
x^{2}-4y\right)  \right)  y^{-1/4}}{\sqrt{\pi}}\int_{2y}^{\infty}\exp\left(
-t\right)  \left(  \frac{\left(  x^{2}+2t\right)  \left(  x^{2}-4y+2t\right)
}{t\left(  t-2y\right)  }\right)  ^{(1+2\nu)/8}\\
& \hspace{1cm}\times\left(  t\left(  x^{2}-4y+2t\right)  \right)
^{-1/4}\left(  P_{-\left(  3+2\nu\right)  /4}^{\left(  1+2\nu\right)
/4}\left(  -\sqrt{\frac{2x^{2}y}{t\left(  x^{2}-4y+2t\right)  }}\right)
\right.  \\
& \hspace{1cm}\left.  -P_{-\left(  3+2\nu\right)  /4}^{\left(  1+2\nu\right)
/4}\left(  \sqrt{\frac{2x^{2}y}{t\left(  x^{2}-4y+2t\right)  }}\right)
\right)  dt\\
& \hspace{1.5cm}\left[  \operatorname{Re}\nu<0,x>0,y>0\right]  .
\end{align*}
\noindent The final example focuses on a limiting case for $D_{\nu}\left(
-x\right)  D_{\mu}\left(  y\right)  $ for which identity (\ref{legendre2}) can
be specialized in terms of the incomplete beta function. Setting $\mu$ at $-1$
in the integral representation (\ref{representation2}) gives an integral
representation for $D_{\nu}\left(  -x\right)  $erfc$\left(  y\right)  $. The
first two Gaussian hypergeometric functions reduce into the incomplete beta
function via the identities (\ref{linear1}), (\ref{linear4}) and
(\ref{incomplete}). The third Gaussian hypergeometric function can also be
simplified into the incomplete beta function given that (\ref{legendre2}) can
be written as%
\begin{align*}
& F\left(  1,b;\tfrac{3}{2};x\right)  =2^{2b-3}x^{-1/2}\left(  1-x\right)
^{\left(  1-2b\right)  /2}\\
& \hspace{1cm}\left(  B\left(  \left(  1+\sqrt{x}\right)
/2;b-1/2,b-1/2\right)  -B\left(  \left(  1-\sqrt{x}\right)
/2;b-1/2,b-1/2\right)  \right)  .
\end{align*}
\noindent This simplification is possible in view of the fact that the degree
and order in the associated Legendre functions in (\ref{legendre2}) both are
at $1/2-b$. The associated Legendre function of the first kind is defined as%
\[
P_{\nu}^{\mu}\left(  x\right)  =\dfrac{1}{\Gamma\left(  1-\mu\right)  }\left(
\dfrac{1+x}{1-x}\right)  ^{\mu/2}F\left(  -\nu,\nu+1;1-\mu;\dfrac{1-x}%
{2}\right)  \text{ for }-1<x<1,
\]
\noindent see (6) on p. 143 of \cite{emoth153}. For identical degree and order
at, for instance, $\nu$ this gives%
\begin{align*}
P_{\nu}^{\nu}\left(  x\right)    & =\dfrac{1}{\Gamma\left(  1-\nu\right)
}\left(  \dfrac{1+x}{1-x}\right)  ^{\nu/2}F\left(  -\nu,\nu+1;1-\nu
;\dfrac{1-x}{2}\right)  \text{ for }-1<x<1,\\
P_{\nu}^{\nu}\left(  x\right)    & =-\dfrac{\nu2^{-\nu}}{\Gamma\left(
1-\nu\right)  }\left(  1-x^{2}\right)  ^{\nu/2}B\left(  \dfrac{1-x}{2}%
;-\nu,-\nu\right)  \text{ for }-1<x<1,
\end{align*}
\noindent on account of the identities (\ref{linear1}) and (\ref{incomplete}).

Using the latter result then implies that all integrands in this integral
representation for $D_{\nu}\left(  -x\right)  $erfc$\left(  y\right)  $ can be
expressed in terms of the incomplete beta function%
\begin{align*}
& D_{\nu}\left(  -x\right)  \text{erfc}\left(  y\right)  =-\dfrac{\nu
\exp\left(  -\tfrac{1}{4}\left(  x^{2}+4y^{2}\right)  \right)  }{2\sqrt{\pi
}\Gamma\left(  1-\nu/2\right)  }\int_{0}^{\infty}\exp\left(  -t\right)
\left(  x^{2}+2y^{2}+2t\right)  ^{\nu/2}\\
& \hspace{1cm}\times\left(  y^{2}+t\right)  ^{-\left(  1+\nu\right)
/2}B\left(  \dfrac{t\left(  x^{2}+2y^{2}+2t\right)  }{\left(  x^{2}+2t\right)
\left(  y^{2}+t\right)  };-\tfrac{\nu}{2},\tfrac{1}{2}\right)  dt\\
& +\dfrac{\left(  1+\nu\right)  \exp\left(  \tfrac{1}{4}\left(  x^{2}%
-4y^{2}\right)  \right)  }{\sqrt{2}\Gamma\left(  -\nu\right)  \Gamma\left(
\left(  3+\nu\right)  /2\right)  }\int_{0}^{\tfrac{1}{2}x^{2}}\exp\left(
-t\right)  \left(  2y^{2}-x^{2}+2t\right)  ^{-\left(  1+\nu\right)  /2}\left(
y^{2}+t\right)  ^{\nu/2}\\
& \hspace{1cm}\times B\left(  \dfrac{t\left(  2y^{2}-x^{2}+2t\right)  }%
{x^{2}y^{2}};\tfrac{1+\nu}{2},-\tfrac{\nu}{2}\right)  dt\\
& +\dfrac{\sqrt{2}\exp\left(  \tfrac{1}{4}\left(  x^{2}-4y^{2}\right)
\right)  }{\sqrt{\pi}\Gamma\left(  -\nu/2\right)  }\int_{\tfrac{1}{2}x^{2}%
}^{\infty}\exp\left(  -t\right)  \left(  2t-x^{2}\right)  ^{-\nu/2}\left(
2y^{2}-x^{2}+2t\right)  ^{-\left(  1+\nu\right)  /2}\\
& \hspace{1cm}\times\left(  t\left(  2y^{2}-x^{2}+2t\right)  -x^{2}%
y^{2}\right)  ^{\nu/2}\left\{  B\left(  \dfrac{1}{2}+\dfrac{1}{2}\sqrt
{\dfrac{x^{2}y^{2}}{t\left(  2y^{2}-x^{2}+2t\right)  }};-\tfrac{\nu}%
{2},-\tfrac{\nu}{2}\right)  \right. \\
& \hspace{1cm}\left.  -B\left(  \dfrac{1}{2}-\dfrac{1}{2}\sqrt{\dfrac
{x^{2}y^{2}}{t\left(  2y^{2}-x^{2}+2t\right)  }};-\tfrac{\nu}{2},-\tfrac{\nu
}{2}\right)  \right\}  dt\\
& \hspace{1.5cm}\left[  -2<\operatorname{Re}\nu<0,x\geqslant0,y>0\right]  .
\end{align*}

\end{document}